%% file: Pochtisovsem.tex
\documentclass[11pt,reqno]{amsart}
\usepackage{amssymb,amscd,amsbsy,mathrsfs}
\input{classstartscr}

\usepackage{amsmath}

\makeatletter
\def\upintkern@{\mkern-7mu\mathchoice{\mkern-3.5mu}{}{}{}}
\def\upintdots@{\mathchoice{\mkern-4mu\@cdots\mkern-4mu}%
 {{\cdotp}\mkern1.5mu{\cdotp}\mkern1.5mu{\cdotp}}%
 {{\cdotp}\mkern1mu{\cdotp}\mkern1mu{\cdotp}}%
 {{\cdotp}\mkern1mu{\cdotp}\mkern1mu{\cdotp}}}

\newcommand{\UpMultiIntegral}[1]{%
  \edef\ints@c{\noexpand\upintop
    \ifnum#1=\z@\noexpand\upintdots@\else\noexpand\upintkern@\fi
    \ifnum#1>\tw@\noexpand\upintop\noexpand\upintkern@\fi
    \ifnum#1>\thr@@\noexpand\upintop\noexpand\upintkern@\fi
    \noexpand\upintop
    \noexpand\ilimits@
  }%
  \futurelet\@let@token\ints@a
}
\makeatother

\DeclareFontFamily{OMX}{mdbch}{}
\DeclareFontShape{OMX}{mdbch}{m}{n}{ <->s * [0.8]  mdbchr7v }{}
\DeclareFontShape{OMX}{mdbch}{b}{n}{ <->s * [0.8]  mdbchb7v }{}
\DeclareFontShape{OMX}{mdbch}{bx}{n}{<->ssub * mdbch/b/n}{}

\DeclareSymbolFont{uplargesymbols}{OMX}{mdbch}{m}{n}
\SetSymbolFont{uplargesymbols}{bold}{OMX}{mdbch}{b}{n}
\DeclareMathSymbol{\upintop}{\mathop}{uplargesymbols}{82}
\DeclareMathSymbol{\upointop}{\mathop}{uplargesymbols}{"48}

\DeclareFontEncoding{MDB}{}{}
\DeclareFontFamily{MDB}{mdbch}{}
\DeclareFontShape{MDB}{mdbch}{m}{n}{ <->s * [0.8]  mdbchrmb }{}
\DeclareFontShape{MDB}{mdbch}{b}{n}{ <->s * [0.8]  mdbchbmb }{}
\DeclareFontShape{MDB}{mdbch}{bx}{n}{<->ssub * mdbch/b/n}{}
\DeclareFontSubstitution{MDB}{cmr}{m}{n}
\DeclareSymbolFont{mathdesignB}{MDB}{mdbch}{m}{n}%
\SetSymbolFont{mathdesignB}{bold}{MDB}{mdbch}{b}{n}%
\DeclareMathSymbol{\upintclockwise}{\mathop}{mathdesignB}{128}
\DeclareMathSymbol{\upointclockwise}{\mathop}{mathdesignB}{130}
\DeclareMathSymbol{\upointctrclockwise}{\mathop}{mathdesignB}{132}
\DeclareMathSymbol{\upoiint}{\mathop}{mathdesignB}{134}
\DeclareMathSymbol{\upoiiint}{\mathop}{mathdesignB}{136}

\makeatletter
\newcommand{\upint}{\DOTSI\upintop\ilimits@}
\newcommand{\upoint}{\DOTSI\upointop\ilimits@}
\makeatother

\theoremstyle{remark}

\newtheorem*{rem*}{Remark}

\newcommand\fM{\frak M}

\newcommand\dg{\frak D}

\newcommand{\fI}{{\frak I}}

\newcommand\mB{\mathcal{B}}

\begin{document}

\newcommand{\vse}{\vspace{.2in}}
\numberwithin{equation}{section}

\title{Functions of almost commuting operators and an extension of the Helton--Howe trace formula}

\author{A.B. Aleksandrov and V.V. Peller}
\thanks{The first author is partially supported by RFBR grant 14-01-00198;
the second author is partially supported by NSF grant DMS 1300924}
\thanks{Corresponding author: V.V. Peller; email: peller@math.msu.edu}

\begin{abstract} 
Let $A$ and $B$ be almost commuting (i.e., the commutator
$AB-BA$ belongs to trace class) self-adjoint operators. We construct a functional calculus
$\f\mapsto\f(A,B)$ for functions $\f$ in the Besov class $B_{\be,1}^1(\R^2)$. This functional calculus is linear, the operators $\f(A,B)$ and $\psi(A,B)$ almost commute for 
$\f,\,\psi\in B_{\be,1}^1(\R^2)$, and $\f(A,B)=u(A)v(B)$ whenever $\f(s,t)=u(s)v(t)$.   We extend the Helton--Howe trace formula for arbitrary functions in $B_{\be,1}^1(\R^2)$. 
The main tool is triple operator integrals with integrands in Haagerup-like tensor products of $L^\be$ spaces.
\end{abstract}

\maketitle

\tableofcontents

\setcounter{section}{0}
\section{Introduction}  
\label{intro}

\

We are going to construct in this paper a functional calculus 
$\f\mapsto\f(A,B)$ for a pair of almost commuting
self-adjoint operators $A$ and $B$. We would like to define such a functional calculus for a substantial class of functions $\f$ of two real variables. Recall that $A$ and $B$ are called {\it almost commuting} if their {\it commutator} $[A,B]\df AB-BA$ belongs to trace class. Certainly, such a functional calculus cannot be multiplicative unless $A$ and $B$ commute. What we would like to have is that the functions of $A$ and $B$ should almost commute, i.e., $[\f(A,B),\psi(A,B)]\in\bS_1$ for arbitrary $\f$ and $\psi$ in our class of functions.
 
Unless otherwise specified, throughout the paper we deal with bounded operators.

Almost commuting self-adjoint operators appear on many occasions. We briefly dwell on two important cases: hyponormal operators and Toeplitz operators.

We can associate with a pair of self-adjoint operators $A$ and $B$ the operator 
$T=A+{\rm i}B$. It is easy to see that $A$ and $B$ almost commute if and only if the {\it selfcommutator} $[T^*,T]$ belongs to trace class. Recall that an operator $T$ is called {\it hyponormal} if $[T^*,T]\ge\0$. It was discovered by Berger and Shaw in \cite{BeSh}
that if $T=A+iB$ is a multicyclic hyponormal operator, then $A$ and $B$ almost commute. Recall that $T$ is called {\it multicyclic} ($m$-{\it multicyclic}) if there are finitely many vectors $x_1,\cdots,x_m$ such that the linear span of 
$$
\{f(T)x_j:~1\le j\le m,\quad f~\mbox{ is a rational function with poles off the spectrum}~\s(T)\}
$$
is dense. Moreover, it was proved in \cite{BeSh} that if $T$ is an $m$-multicyclic hyponormal operator, then
$$
\big\|\,[T^*,T]\,\big\|_{\bS_1}\le\frac{m}\pi\,\m_2(\s(T)),
$$
where $\m_2$ is planar Lebesgue measure.

Another important example of almost commuting pairs of self-adjoint operators is given by Toeplitz operators. Suppose that $f$ and $g$ are bounded real functions on the unit circle $\T$ such that $f$  belongs to the Besov class $B_p^{1/p}$ and $g$ belongs to
$B_{p'}^{1/p'}$, $1<p<\be$, $1/p+1/p'=1$ (or $f\in L^\be$ and $g\in B_1^1$), then
the Toeplitz operators $T_f$ and $T_g$ form a pair of almost commuting self-adjoint operators. Indeed, in this case it is easy to see that
$$
\big[T_f,T_g\big]=H^*_{\ov g}H_f-H^*_{\ov f}H_g,
$$
where $H_f$, $H_g$, $H_{\ov f}$ and $H_{\ov g}$ are Hankel operators. The fact that 
$[T_f,T_g\big]\in\bS_1$ follows from the description of Hankel operators of Schatten--von Neumann class $\bS_p$ that was obtained in \cite{Pe1}, see also the book \cite{Pe5}, Ch. 6.

The polynomial calculus for a pair $A$ and $B$ of almost commuting self-adjoint operators can be defined in the following way. For a polynomial $\f$ of the form $\f(s,t)=\sum_{j,k}a_{jk}s^jt^k$, the operator $\f(A,B)$ is defined by
$$
\f(A,B)\df\sum_{j,k}a_{jk}A^jB^k.
$$
It can easily be verified that if $\f$ and $\psi$ are polynomials of two real variables then $[\f(A,B),\psi(A,B)]\in\bS_1$.

In \cite{HH} the following trace formula was obtained for bounded almost commuting self-adjoint operators $A$ and $B$:
\bay
\label{HeHo}
\trace\big({\rm i}\big(\f(A,B)\psi(A,B)-\psi(A,B)\f(A,B)\big)\big)=
\int_{\R^2}\left(\frac{\partial\f}{\partial x}\frac{\partial\psi}{\partial y}-
\frac{\partial\f}{\partial y}\frac{\partial\psi}{\partial x}\right)dP,
\ey
where $P$ is a signed Borel compactly supported measure that corresponds to the pair $(A,B)$. The formula holds for polynomials $\f$ and $\psi$. 
It was shown in \cite{Pin2} that the signed measure $P$ is absolutely continuous with respect to planar Lebesgue measure and
$$
dP(x,y)=\frac1{2\pi}g(x,y)\,d\m_2(x,y),
$$
where $g$ is the {\it Pincus principal function}, which was introduced in \cite{Pin1}. We refer the reader to \cite{Cl} for more detailed information.

Note that Helton and Howe extended the polynomial calculus in \cite{HH} to the class of infinitely differentiable functions and they proved that trace formula \rf{HeHo} holds for such functions.

In \cite{CP} the polynomial functional calculus for almost commuting self-adjoint operators was extended to a functional calculus for the class of functions $\f=\F\o$ that are Fourier transforms of complex Borel measures $\o$ on $\R^2$ satisfying
$$
\int_{\R^2}(1+|t|)(1+|s|)\,d|\o|(s,t)<\be,
$$
and the Helton--Howe trace formula \rf{HeHo} was extended to the class of such functions.

The spectral theorem for pairs of commuting self-adjoint operators, associates with such a pair $A$ and $B$ the spectral measure $E_{A,B}$ on Borel subsets of the plain $\R^2$.
This allows one to construct a linear and multiplicative functional calculus 
$$
\f\mapsto\f(A,B)=\int_{\R^2}\f(x,y)\,dE_{A,B}(x,y)
$$
for the class of bounded Borel functions on the plane $\R^2$. The support of the 
spectral measure $E_{A,B}$ coincides with the joint spectrum of the pair $(A,B)$.

If $A$ and $B$ are noncommuting self-adjoint operators, we can define functions of $A$ and $B$ in terms of double operator integrals
\bay
\label{isch}
\f(A,B)\df\iint_{\R^2}\f(x,y)\,dE_A(x)\,dE_B(y).
\ey
However, unlike in the case of commuting self-adjoint operators, we cannot define
functions $\f(A,B)$ for arbitrary bounded Borel functions $\f$. 
Such double operator integrals can be defined
for functions $\f$ that are Schur multipliers with respect to the spectral measure $E_A$ and $E_B$ of the operators $A$ and $B$. The theory of double operator integrals was developed by Birman and Solomyak \cite{BS} (we also refer the reader to \cite{Pe2} and \cite{AP} for double operator integrals and Schur multipliers).

The problem of constructing a rich functional calculus for almost commiting self-adjoint operators, which would extend the functional calculus constructed in \cite{CP} and for which trace formula \rf{HeHo} would still hold was considered in \cite{Pe4}. The problem was to find a big class of functions ${\mathcal C}$ on $\R^2$ and construct a functional calculus $\f\mapsto\f(A,B)$, $\f\in{\mathcal C}$, that has the following properties:

{\it
{\em(i)}~ the functional calculus $\f\mapsto\f(A,B)$, $\f\in{\mathcal C}$, is linear;

{\em(ii)}~ if $\f(s,t)=u(s)v(t)$, then $\f(A,B)=u(A)v(B)$;

{\em(iii)}~ if $\f,\,\psi\in{\mathcal C}$, then  
$\f(A,B)\psi(A,B)-\psi(A,B)\f(A,B)\in\bS_1$;

{\em(iv)}~ formula {\em\rf{HeHo}} holds for arbitrary $\f$ and $\psi$ in ${\mathcal C}$.
}

Note that the right-hand side of \rf{HeHo} makes sense for arbitrary Lipschitz functions $\f$ and $\psi$. However, it was established
in \cite{Pe4}  that a functional calculus satisfying (i) - (iii) cannot be defined for all continuously differentiable functions. This was deduced from the trace class criterion for Hankel operators (see \cite{Pe1} and \cite{Pe5}).

On the other hand, in \cite{Pe4} estimates of \cite{Pe2} and \cite{Pe3}  were used to construct a functional calculus satisfying (i) - (iv)  for the class 
${\mathcal C}=\big(L^\be(\R)\hat\otimes B_{\be,1}^1(\R)\big)\bigcap\big(B_{\be,1}^1(\R)\hat\otimes L^\be(\R)\big)$. Here $\hat\otimes$ stands for projective tensor product and $B_{\be,1}^1(\R)$ is a Besov class (see \S~\ref{Pre} for a brief introduction to Besov classes).

In this paper we considerably enlarge the class $\big(L^\be(\R)\hat\otimes B_{\be,1}^1(\R)\big)\!\bigcap\!\big(B_{\be,1}^1(\R)\hat\otimes L^\be(\R)\big)$ and construct a functional calculus satisfying (i) - (iv) for the Besov class
$B_{\be,1}^1(\R^2)$ of functions of two variables.

It was observed in \cite{ANP} that the Besov space $B_{\infty,1}^1(\R^2)$ 
of functions on $\R^2$ is contained in the space of Schur multipliers with respect to compactly supported spectral measures on $\R$, and so for 
$\f\in B_{\infty,1}^1(\R^2)$, the operator $\f(A,B)$ is well defined by \rf{isch} for bounded self-adjoint operators $A$ and $B$.

The results of this paper were announced in \cite{AP2}.

In \S~\ref{tri} we deal with triple operator integrals. We consider triple operator integrals with integrands in the Haagerup tensor product of $L^\be$ spaces. It turns out that for our purpose such triple operator integrals cannot be used. We define in 
\S~\ref{tri} Haagerup-like tensor products of the first kind and of the second kind. Then we define triple operator integrals with symbols in such Haagerup-like tensor products.

We use such triple operator integrals in \S~\ref{comf} to obtain a representation of 
commutators $[\f(A,B),\psi(A,B)]$, $\f,\,\psi\in B_{\be,1}^1(\R^2)$ in terms of triple operator integrals. This allows us to estimate trace norms of such commutators.

In \S~\ref{HelHow} we use the results of \S~\ref{comf} to obtain an extension of the Helton--Howe trace formula for functions in $B_{\be,1}^1(\R^2)$.

In the final section we state open problems.

Finally, we give in \S~\ref{Pre} brief introductions to Besov spaces and double operator integrals.

\

\section{Preliminaries}
\label{Pre}

\

In this section we collect necessary information on Besov spaces and double operator integrals. 

\medskip

{\bf 2.1. Besov classes of functions on Euclidean spaces.} The technique of Little\-wood--Paley type expansions of functions or distributions on Euclidean spaces 
is a very important tool in Harmonic Analysis. 

Let $w$ be an infinitely differentiable function on $\R$ such
that
\bay
\label{w}
w\ge0,\quad\supp w\subset\left[\frac12,2\right],\quad\mbox{and} \quad w(s)=1-w\left(\frac s2\right)\quad\mbox{for}\quad s\in[1,2].
\ey

We define the functions $W_n$, $n\in\Z$, on $\R^d$ by 
$$
\big(\F W_n\big)(x)=w\left(\frac{\|x\|_2}{2^n}\right),\quad n\in\Z, \quad x=(x_1,\cdots,x_d),
\quad\|x\|_2\df\left(\sum_{j=1}^dx_j^2\right)^{1/2},
$$
where $\F$ is the {\it Fourier transform} defined on $L^1\big(\R^d\big)$ by
$$
\big(\F f\big)(t)=\!\int\limits_{\R^d} f(x)e^{-{\rm i}(x,t)}\,dx,\!\quad 
x=(x_1,\cdots,x_d),
\quad t=(t_1,\cdots,t_d), \!\quad(x,t)\df \sum_{j=1}^dx_jt_j.
$$
Clearly,
$$
\sum_{n\in\Z}(\F W_n)(t)=1,\quad t\in\R^d\setminus\{0\}.
$$

With each tempered distribution $f\in{\mathscr S}^\prime\big(\R^d\big)$, we
associate the sequence $\{f_n\}_{n\in\Z}$,
\bay
\label{fn}
f_n\df f*W_n.
\ey
The formal series
$
\sum_{n\in\Z}f_n
$
is a Littlewood--Paley type expansion of $f$. This series does not necessarily converge to $f$. 

Initially we define the (homogeneous) Besov class $\dot B^s_{p,q}\big(\R^d\big)$,
$s>0$, $1\le p,\,q\le\be$, as the space of all
$f\in{\mathscr S}^\prime(\R^n)$
such that
\bay
\label{Wn}
\{2^{ns}\|f_n\|_{L^p}\}_{n\in\Z}\in\ell^q(\Z)
\ey
and put
$$
\|f\|_{B^s_{p,q}}\df\big\|\{2^{ns}\|f_n\|_{L^p}\}_{n\in\Z}\big\|_{\ell^q(\Z)}.
$$
According to this definition, the space $\dot B^s_{p,q}(\R^n)$ contains all polynomials
and all polynomials $f$ satisfy the equality $\|f\|_{B^s_{p,q}}=0$. Moreover, the distribution $f$ is determined by the sequence $\{f_n\}_{n\in\Z}$
uniquely up to a polynomial. It is easy to see that the series 
$\sum_{n\ge0}f_n$ converges in ${\mathscr S}^\prime(\R^d)$.
However, the series $\sum_{n<0}f_n$ can diverge in general. It can easily be proved that the series
\bay
\label{ryad}
\sum_{n<0}\frac{\partial^r f_n}{\partial x_1^{r_1}\cdots\partial x_d^{r_d}},\qquad \mbox{where}\quad r_j\ge0,\quad\mbox{for}\quad
1\le j\le d,\quad\sum_{j=1}^dr_j=r,
\ey
converges uniformly on $\R^d$ for every nonnegative integer
$r>s-d/p$. Note that in the case $q=1$ the series \rf{ryad}
converges uniformly, whenever $r\ge s-d/p$.

Now we can define the modified (homogeneous) Besov class $B^s_{p,q}\big(\R^d\big)$. We say that a distribution $f$
belongs to $B^s_{p,q}(\R^d)$ if \rf{Wn} holds and
$$
\frac{\partial^r f}{\partial x_1^{r_1}\cdots\partial x_d^{r_d}}
=\sum_{n\in\Z}\frac{\partial^r f_n}{\partial x_1^{r_1}\cdots\partial x_d^{r_d}},\quad
\mbox{whenever}\quad 
r_j\ge0,\quad\mbox{for}\quad
1\le j\le d,\quad\sum_{j=1}^dr_j=r.
$$
in the space ${\mathscr S}^\prime\big(\R^d\big)$, where $r$ is
the minimal nonnegative integer such that $r>s-d/p$ ($r\ge s-d/p$ if $q=1$). Now the function $f$ is determined uniquely by the sequence $\{f_n\}_{n\in\Z}$ up
to a polynomial of degree less than $r$, and a polynomial $g$ belongs to 
$B^s_{p,q}\big(\R^d\big)$
if and only if $\deg g<r$.

In this paper we deal with Besov classes 
$B_{\be,1}^1(\R^d)$. They can also be defined in the following way:

Let $X$ be the set of all continuous functions $f\in L^\be(\R^d)$ such that $|f|\le1$ and 
$\supp\F f\subset\{\xi\in\R^d:~\|\xi\|\le1\}$.
Then 
$$
B^1_{\be1}(\R^d)=\left\{c+\sum_{n=1}^\be \a_n\s_n^{-1}(f_n(\s_nx)-f(0)):
~c\in\C,~f_n\in X, ~\s_n>0,~\sum_{n=1}^\be|\a_n|<\be\right\}.
$$

Note that the functions $f_\s$, $f_\s(x)=f(\s x)$, $x\in\R^d$, have the following properties: \lb$f_\s\in L^\be(\R^d)$ and 
$\supp\F f\subset\{\xi\in\R^d:~\|\xi\|\le\s\}$. Such functions can be characterized by the following Paley--Wiener--Schwartz type theorem  (see \cite{R}, Theorem 7.23 and exercise 15 of Chapter 7):

{\it Let $f$ be a continuous function
on $\R^d$ and let $M,\,\s>0$. The following statements are equivalent:

{\em(i)} $|f|\le M$ and $\supp\F f\subset\{\xi\in\R^d:\|\xi\|\le\s\}$;

{\em(ii)} $f$ is a restriction to $\R^d$ of an entire function on $\C^d$ such that 
$$
|f(z)|\le Me^{\s\|\im z\|}
$$
for all $z\in\C^d$.}

We refer the reader to \cite{Pee} and \cite{Tr} for more detailed information on Besov spaces.

\medskip

{\bf 2.2. Besov classes of periodic functions.} Studying periodic functions on $\R^d$ is equivalent to studying functions on the $d$-dimensional torus $\T^d$. To define Besov spaces on $\T^d$, we consider a function $w$ satisfying \rf{w} and define the trigonometric polynomials $W_n$, $n\ge0$, by
$$
W_n(\z)\df\sum_{j\in\Z^d}w\left(\frac{|j|}{2^n}\right)\z^j,\quad n\ge1,
\quad W_0(\z)\df\sum_{\{j:|j|\le1\}}\z^j,
$$
where 
$$
\z=(\z_1,\cdots,\z_d)\in\T^d,\quad j=(j_1,\cdots,j_d),\quad\mbox{and}\quad
|j|=\big(|j_1|^2+\cdots+|j_d|^2\big)^{1/2}.
$$
For a distribution $f$ on $\T^d$ we put
$$
f_n=f*W_n,\quad n\ge0,
$$
and we say that $f$ belongs the Besov class $B_{p,q}^s(\T^d)$, $s>0$, 
$1\le p,\,q\le\be$, if
\bay
\label{Bperf}
\big\{2^ns\|f_n\|_{L^p}\big\}_{n\ge0}\in\ell^q.
\ey

Note that locally the Besov space $B_{p,q}^s(\R^d)$ coincides with the Besov space
$B_{p,q}^s$ of periodic functions on $\R^d$.

\medskip

{\bf 2.3. Double operator integrals.}
In this subsection we give a brief introduction to double  operator integrals. Double operator integrals appeared in the paper \cite{DK} by Daletskii and S.G. Krein. Later the beautiful theory of double operator integrals was developed by Birman and Solomyak in \cite{BS}, \cite{BS2}, and \cite{BS3}.

Let $(\X,E_1)$ and $(\Y,E_2)$ be spaces with spectral measures $E_1$ and $E_2$
on a Hilbert space $\h$. The idea of Birman and Solomyak is to define first
double operator integrals
\bay
\label{doi}
\int\limits_\X\int\limits_\Y\Phi(x,y)\,d E_1(x)T\,dE_2(y),
\ey
for bounded measurable functions $\Phi$ and operators $T$
of Hilbert Schmidt class $\bS_2$. Consider the spectral measure $\E$ whose values are orthogonal projections on the Hilbert space $\bS_2$, which is defined by
$$
\E(\L\times\D)T=E_1(\L)TE_2(\D),\quad T\in\bS_2,
$$
$\L$ and $\D$ being measurable subsets of $\X$ and $\Y$. It was shown in \cite{BS4} that $\E$ extends to a spectral measure on
$\X\times\Y$. If $\Phi$ is a bounded measurable function on $\X\times\Y$, we define the double operator integral \rf{doi} by
$$
\int\limits_\X\int\limits_\Y\Phi(x,y)\,d E_1(x)T\,dE_2(y)\df
\left(\,\,\int\limits_{\X\times\Y}\Phi\,d\E\right)T.
$$
Clearly,
$$
\left\|\int\limits_\X\int\limits_\Y\Phi(x,y)\,dE_1(x)T\,dE_2(y)\right\|_{\bS_2}
\le\|\Phi\|_{L^\be}\|T\|_{\bS_2}.
$$
If
$$
\int\limits_\X\int\limits_\Y\Phi(x,y)\,d E_1(x)T\,dE_2(y)\in\bS_1
$$
for every $T\in\bS_1$, we say that $\Phi$ is a {\it Schur multiplier of $\bS_1$ associated with
the spectral measures $E_1$ and $E_2$}.

To define double operator integrals of the form \rf{doi} for bounded linear operators $T$,
we consider the transformer
$$
Q\mapsto\int\limits_{\Y}\int\limits_{\X}\Phi(y,x)\,d E_2(y)\,Q\,dE_1(x),\quad Q\in\bS_1,
$$
and assume that the function $(y,x)\mapsto\Phi(y,x)$ is a Schur multiplier of $\bS_1$ with associated with $E_2$ and $E_1$.

In this case the transformer
\bay
\label{tra}
T\mapsto\int\limits_\X\int\limits_\Y\Phi(x,y)\,d E_1(x)T\,dE_2(y),\quad T\in \bS_2,
\ey
extends by duality to a bounded linear transformer on the space of bounded linear operators on $\h$
and we say that the function $\Phi$ is {\it a Schur multiplier (with respect to $E_1$ and $E_2$) of the space of bounded linear operators}.
We denote the space of such Schur multipliers by $\fM(E_1,E_2)$.
The norm of $\Phi$ in $\fM(E_1,E_2)$ is, by definition, the norm of the
transformer \rf{tra} on the space of bounded linear operators.


It was observed in \cite{BS5} that if $A$ and $B$ are self-adjoint operators and if $f$ is a continuously differentiable
function on $\R$ such that the divided difference $\dg f$,
$$
\big(\dg f\big)(x,y)=\frac{f(x)-f(y)}{x-y},
$$
is a Schur multiplier
with respect to the spectral measures of $A$ and $B$, then for an arbitrary bounded line operator $Q$ the following formula holds
$$
f(A)Q-Qf(B)=\iint\big(\dg f\big)(x,y)\,dE_{A}(x)(AQ-QB)\,dE_B(y)
$$
and
$$
\|f(A)Q-Qf(B)\|\le\const\|\dg f\|_{\fM(E_A,E_{B})}\|AQ-QB\|.
$$
The same inequality holds if we replace the operator norm with a norm in a separable symmetrically normed ideal (see \cite{GK}), in particular, in the Schatten--von Neumann norms $\bS_p$.

It was established in \cite{Pe2} (see also \cite{Pe3}) that if $f$ belongs to the Besov class $B_{\be,1}^1(\R)$, then the divided difference $\dg f\in \fM(E_1,E_2)$
for arbitrary Borel spectral $E_1$ and $E_2$, and so 
\bay
\label{BesSch}
\|f(A)Q-Qf(B)\|_\fI\le\const\|f\|_{B_{\be,1}^1}\|AQ-QB\|_\fI
\ey
for arbitrary self-adjoint operators $A$ and $B$ and an arbitrary separable symmetrically normed ideal $\fI$.

There are different characterizations of the space $\fM(E_1,E_2)$ of Schur multipliers, see \cite{Pe2} and \cite{Pis}. In particular, a function $\Phi$ is a Schur multiplier if and only if it belongs to the {\it Haagerup tensor product} $L^\be(E_1)\!\otimes_{\rm h}\!L^\be(E_2)$, which is, by definition, 
the space of functions $\Phi$ of the form 
\bay
\label{FiH}
\Phi(x,y)=\sum_{j\ge0}\f_j(x)\psi_j(y),
\ey
where $\f_j\in L^\be(E_1)$, $\psi_j\in L^\be(E_2)$ and
$$
\{\f_j\}_{j\ge0}\in L_{E_1}^\be(\ell^2)\quad\mbox{and}\quad
\{\psi_j\}_{j\ge0}\in L_{E_2}^\be(\ell^2).
$$
The {\it norm of $\Phi$ in $L^\be(E_1)\!\otimes_{\rm h}\!L^\be(E_2)$ is defined as} the infimum of
$$
\big\|\{\f_j\}_{j\ge0}\big\|_{L_{E_1}^\be(\ell^2)}
\big\|\{\psi_j\}_{j\ge0}\big\|_{L_{E_2}^\be(\ell^2)}
$$
over all representations of $\Phi$ of the form \rf{FiH}. Here
$$
\big\|\{\f_j\}_{j\ge0}\big\|_{L_{E_1}^\be(\ell^2)}\df
\Big\|\sum_{j\ge0}|\f_j|^2\Big\|_{L^\be(E_1)}^{1/2}\quad\!\!\mbox{and}\quad\!\!
\big\|\{\psi_j\}_{j\ge0}\big\|_{L_{E_1}^\be(\ell^2)}\df
\Big\|\sum_{j\ge0}|\psi_j|^2\Big\|_{L^\be(E_2)}^{1/2}.
$$
It can easily be verified that if $\Phi\in L^\be(E_1)\!\otimes_{\rm h}\!L^\be(E_2)$, then $\Phi\in\fM(E_1,E_2)$ and
\bay
\label{skh}
\iint\Phi(x,y)\,dE_1(x)T\,dE_2(y)=
\sum_{j\ge0}\Big(\int\f_j\,dE_1\Big)T\Big(\int\psi_j\,dE_2\Big)
\ey
and the series on the right converges in the weak operator topology.
It is also easy to see that the series on the right converges in the weak operator topology and
$$
\|\Phi\|_{\fM(E_1,E_2)}\le\|\Phi\|_{L^\be(E_1)\otimes_{\rm h}L^\be(E_2)}.
$$


Let us also mention the following sufficient condition:

\medskip

{\it If a function $\Phi$ on $\X\times\Y$ belongs to the {\it projective tensor
product}
$L^\be(E_1)\hat\otimes L^\be(E_2)$ of $L^\be(E_1)$ and $L^\be(E_2)$ (i.e., $\Phi$ admits a representation of the form {\em\rf{FiH}}
whith \lb$\f_j\in L^\be(E_1)$, $\psi_j\in L^\be(E_2)$, and
$$
\sum_{j\ge0}\|\f_j\|_{L^\be}\|\psi_j\|_{L^\be}<\be),
$$
then $\Phi\in\fM(E_1,E_2)$ and}
\bay
\label{dous}
\|\Phi\|_{\fM(E_1,E_2)}\le\sum_{j\ge0}\|\f_j\|_{L^\be}\|\psi_j\|_{L^\be}.
\ey

For such functions $\Phi$, formula \rf{skh} holds
and the series on the right-hand side of \rf{skh} converges absolutely in the norm.

\medskip

{\bf2.4. Functions of noncommuting self-adjoint operators.} Let $A$ and $B$ be self-adjoint operators on Hilbert space and let $E_A$ and $E_B$ be their spectral measures. Suppose that $f$ is a function of two variables that is defined at least on $\s(A)\times\s(B)$. As we have already mentioned in the introduction, if
 $f$ is a Schur multiplier with respect to the pair $(E_A,E_B)$, we define the function $f(A,B)$ of $A$ and $B$ by
\bay
\label{fAB}
f(A,B)\df\iint f(x,y)\,dE_A(x)\,dE_B(y).
\ey
Note that this functional calculus $f\mapsto f(A,B)$ is linear, but not multiplicative.

If we consider functions of bounded operators, without loss of generality we may deal with periodic functions with a sufficiently large period. Clearly, we can rescale the problem and assume that our functions are $2\pi$-periodic in each variable.

If $f$ is a trigonometric polynomial of degree $N$, we can represent $f$ in the form
$$
f(x,y)=\sum_{j=-N}^Ne^{{\rm i}jx}\left(\sum_{k=-N}^N\hat f(j,k)e^{{\rm i}ky}\right).
$$
Thus 
$f$ belongs to the projective tensor product $L^\be\hat\otimes L^\be$ and
$$
\|f\|_{L^\be\hat\otimes L^\be}\le\sum_{j=-N}^N\sup_y
\left|\sum_{k=-N}^N\hat f(j,k)e^{{\rm i}ky}\right|
\le(1+2N)\|f\|_{L^\be}
$$
It follows easily from \rf{Bperf} that every periodic function $f$ of Besov class $B_{\be1}^1$ of periodic functions belongs to 
$L^\be\hat\otimes L^\be$, and so the operator $f(A,B)$ is well defined by \rf{fAB}.

\

\section{Triple operator integrals}
\label{tri}

\

Multiple operator integrals were considered by several mathematicians, see \cite{Pa}, \cite{St}. However, those definitions required very strong restrictions on the classes of functions that can be integrated. In \cite{Pe6} multiple operator integrals were defined for functions that belong to the (integral) projective tensor product of $L^\be$ spaces. Later in \cite{JTT} multiple operator integrals were defined for Haagerup tensor products of $L^\be$ spaces.

In this paper we deal with triple operator integrals. 
We consider here both approaches given in \cite{Pe6} and \cite{JTT}.

It turns out, however that none of these approaches helps in our situation. That is why we define Haagerup-like tensor products of the first kind and of the second kind and define triple operator integrals whose integrands belong to such Haagerup-like tensor products.

Let $E_1$, $E_2$, and $E_3$ be spectral measures on Hilbert space and let $T$ and $R$ be bounded linear operators on Hilbert space. Triple operator integrals are expressions of the following form:
\bay
\label{troi}
\int\limits_{\X_1}\int\limits_{\X_2}\int\limits_{\X_3} 
\Psi(x_1,x_2,x_3)\,dE_1(x_1)T\,dE_2(x_2)R\,dE_3(x_3).
\ey
Such integrals make sense under certain assumptions on $\Psi$, $T$, and $R$. The function $\Psi$ will be called the {\it integrand} of the triple operator integral.

Recall that the {\it projective tensor product} 
$L^\be(E_1)\hat\otimes L^\be(E_2)\hat\otimes L^\be(E_3)$ can be defined as the class of function $\Psi$ of the form
\bay
\label{pred}
\Psi(x_1,x_2,x_3)=\sum_n\f_n(x_1)\psi_n(x_2)\chi_n(x_3)
\ey
such that
\bay
\label{norma}
\sum_n\|\f_n\|_{L^\be(E_1)}\|\psi_n\|_{L^\be(E_2)}\|\chi_n\|_{L^\be(E_3)}<\be.
\ey
The norm $\|\Psi\|_{L^\be\hat\otimes L^\be\hat\otimes L^\be}$ of $\Psi$ is, by definition, the infimum of the left-hand side of \rf{norma} over all representations of the form \rf{pred}.

For $\Psi\in L^\be(E_1)\hat\otimes L^\be(E_2)\hat\otimes L^\be(E_3)$ of the form 
\rf{pred} the triple operator integral \rf{troi} was defined in \cite{Pe4} by
\begin{align}
\label{otoi}
\iiint \Psi(x_1,x_2,x_3)&\,dE_1(x_1)T\,dE_2(x_2)R\,dE_3(x_3)\nonumber\\
=&\sum_n\left(\int\f_n\,dE_1\right)T\left(\int\psi_n\,dE_2\right)R\left(\int\chi_n\,dE_3\right).
\end{align}
Clearly, \rf{norma} implies that the series on the right converges absolutely in the norm. The right-hand side of \rf{otoi} does not depend on the choice of a representation of the form \rf{pred}. Clearly,
$$
\left\|\iiint \Psi(x_1,x_2,x_3)\,dE_1(x_1)T\,dE_2(x_2)R\,dE_3(x_3)\right\|
\le\|\Psi\|_{L^\be\hat\otimes L^\be\hat\otimes L^\be}\|T\|\cdot\|R\|.
$$
Note that for $\Psi\in L^\be(E_1)\hat\otimes L^\be(E_2)\hat\otimes L^\be(E_3)$, triple operator integrals have the following properties:
\bay
\label{bep}
T\in\mB(\h),\quad R\in\bS_p,\quad1\le p<\be,\quad\Longrightarrow\quad 
\iiint\Psi\,dE_1T\,dE_2R\,dE_3\in\bS_p
\ey
and
\bay
\label{pq}
T\in\bS_p,\!\quad\! R\in\bS_q,\quad \frac1p+\frac1q\le1\!
\!\quad\Longrightarrow\!\quad 
\iiint\!\Psi dE_1TdE_2RdE_3\in\bS_r,\quad\!\frac1r=\frac1p+\frac1q.
\ey

Let us also mention that multiple operator integrals were defined in \cite{Pe6} for 
functions $\Psi$ that belong to the so-called {\it integral projective tensor product} of the corresponding $L^\be$ spaces, which contains the projective tenser product. We refer the reader to \cite{Pe6} for more detail.

We proceed now to the approach to multiple operator integrals based on the Haagerup tensor product of $L^\be$ spaces. We refer the reader to the book \cite{Pis} for detailed \lb information about Haagerup tensor products.
We define the {\it Haagerup tensor product} \lb
$L^\be(E_1)\!\otimes_{\rm h}\!L^\be(E_2)\!\otimes_{\rm h}\!L^\be(E_3)$ as the space of function $\Psi$ of the form
\bay
\label{htr}
\Psi(x_1,x_2,x_3)=\sum_{j,k\ge0}\a_j(x_1)\b_{jk}(x_2)\g_k(x_3),
\ey
where $\a_j$, $\b_{jk}$, and $\g_k$ are measurable functions such that
\bay
\label{ogr}
\{\a_j\}_{j\ge0}\in L_{E_1}^\be(\ell^2), \quad 
\{\b_{jk}\}_{j,k\ge0}\in L_{E_2}^\be({\mathcal B}),\quad\mbox{and}\quad
\{\g_k\}_{k\ge0}\in L_{E_3}^\be(\ell^2),
\ey
where ${\mathcal B}$ is the space of matrices that induce bounded linear operators on $\ell^2$ and this space is equipped with the operator norm. In other words,
$$
\|\{\a_j\}_{j\ge0}\|_{L^\be(\ell^2)}\df
E_1\mbox{-}\ess\sup\left(\sum_{j\ge0}|\a_j(x_1)|^2\right)^{1/2}<\be,
$$
$$
\|\{\b_{jk}\}_{j,k\ge0}\|_{L^\be({\mathcal B})}\df
E_2\mbox{-}\ess\sup\|\{\b_{jk}(x_2)\}_{j,k\ge0}\|_{{\mathcal B}}<\be,
$$
and
$$
\|\{\g_k\}_{k\ge0}\|_{L^\be(\ell^2)}\df
E_3\mbox{-}\ess\sup\left(\sum_{k\ge0}|\g_k(x_3)|^2\right)^{1/2}<\be.
$$
By the sum on the right-hand of \rf{htr} we mean 
$$
\lim_{M,N\to\be}~\sum_{j=0}^N\sum_{k=0}^M\a_j(x_1)\b_{jk}(x_2)\g_k(x_3).
$$
Clearly, the limit exists.

{\it Throughout the paper by $\sum_{j,k\ge0}$, we mean 
$\lim_{M,N\to\be}~\sum_{j=0}^N\sum^M_{k=0}$}.

The norm of $\Psi$ in 
$L^\be\!\otimes_{\rm h}\!L^\be\!\otimes_{\rm h}\!L^\be$ is, by definition, 
the infimum of
$$
\|\{\a_j\}_{j\ge0}\|_{L^\be(\ell^2)}\|\{\b_{jk}\}_{j,k\ge0}\|_{L^\be({\mathcal B})}
\|\{\g_k\}_{k\ge0}\|_{L^\be(\ell^2)}
$$
over all representations of $\Psi$ of the form \rf{htr}.

It is easy to verify that $L^\be\hat\otimes L^\be\hat\otimes L^\be\subset
L^\be\!\otimes_{\rm h}\!L^\be\!\otimes_{\rm h}\!L^\be$.

In \cite{JTT} multiple operator integrals were defined for functions in the Haagerup tensor product of $L^\be$ spaces. Let 
$\Psi\in L^\be\!\otimes_{\rm h}\!L^\be\!\otimes_{\rm h}\!L^\be$ and  
suppose that \rf{htr} and \rf{ogr} hold. The triple operator integral \rf{troi} is defined by
\begin{align}
\label{htraz}
\iiint\Psi(x_1,x_2,x_3)&\,dE_1(x_1)T\,dE_2(x_2)R\,dE_3(x_3)\nonumber\\[.2cm]
=&
\sum_{j,k\ge0}\left(\int\a_j\,dE_1\right)T\left(\int\b_{jk}\,dE_2\right)
R\left(\int\g_k\,dE_3\right)\nonumber\\[.2cm]
=&\lim_{M,N\to\be}~\sum_{j=0}^N\sum_{k=0}^M
\left(\int\a_j\,dE_1\right)T\left(\int\b_{jk}\,dE_2\right)
R\left(\int\g_k\,dE_3\right).
\end{align}

Then (see \cite{JTT} and \cite{ANP2}) series in \rf{htraz} converges in the weak operator topology, the sum of the series does not depend on the choice of a representation, and the following inequality holds:
\bay
\label{opno}
\!\!\!\!\left\|\iiint\Psi(x_1,x_2,x_3)\,dE_1(x_1)T\,dE_2(x_2)R\,dE_3(x_3)\right\|
\le\|\Psi\|_{L^\be\!\otimes_{\rm h}\!L^\be\!\otimes_{\rm h}\!L^\be}
\|T\|\cdot\|R\|.
\ey

However, it was shown in \cite{ANP} that unlike in the case 
$\Phi\in L^\be(E_1)\hat\otimes_{\rm i}L^\be(E_2)\hat\otimes_{\rm i}L^\be(E_3)$, the condition 
that $\Psi$ belongs to $L^\be(E_1)\otimes_{\rm h}\!L^\be(E_2)\otimes_{\rm h}\!L^\be(E_3)$ does not guarantee that if one of the operators $T$ and $R$ is of trace class, then the triple operator integral \rf{troi} belongs to $\bS_1$. The same can be said if we replace trace class $\bS_1$ with the Schatten--von Neumann class $\bS_p$ with $p<2$, see \cite{ANP3} and \cite{ANP2}.

However, in this paper we need estimates of triple operator integrals in the norm of $\bS_1$. Also, we are going to use triple operator integrals whose integrands are the divide differences of functions in $B_{\be,1}^1(\R^2)$ in each variable (see \S~\ref{comf}). It was established in \cite{ANP} and \cite{ANP3} that such divided differences do not have to belong to the Haagerup tensor product of $L^\be$ spaces. That is why we need to modify the notion of the Haagerup tensor product.

In \cite{ANP} and \cite{AP2} the following Haagerup-like tensor products were introduced:

\medskip

{\bf Definition 1.}
{\it A function $\Psi$ is said to belong to the 
Haagerup-like tensor product 
$L^\be(E_1)\otimes_{\rm h}\!L^\be(E_2)\otimes^{\rm h}\!L^\be(E_3)$ 
of the first kind if it admits a representation
\bay
\label{yaH}
\Psi(x_1,x_2,x_3)=\sum_{j,k\ge0}\a_j(x_1)\b_{k}(x_2)\g_{jk}(x_3)
\ey
with $\{\a_j\}_{j\ge0},~\{\b_k\}_{k\ge0}\in L^\be(\ell^2)$ and 
$\{\g_{jk}\}_{j,k\ge0}\in L^\be(\mB)$, where $\mB$ is the space of bounded operators on $\ell^2$}.

\medskip

For a bounded linear operator $R$ and 
for a trace class operator $T$,  the triple operator integral
$$
W=\iint\!\!\upint\Psi(x_1,x_2,x_3)\,dE_1(x_1)T\,dE_2(x_2)R\,dE_3(x_3)
$$
was defined in \cite{ANP}
as the following continuous linear functional on the class of compact operators:
\bay
\label{fko}
Q\mapsto
\trace\left(\left(
\iiint\Psi(x_1,x_2,x_3)\,dE_2(x_2)R\,dE_3(x_3)Q\,dE_1(x_1)
\right)T\right)
\ey

The fact that the linear functional \rf{fko} is continuous on the class of compact operators is a consequence of inequality \rf{opno}, which also implies the following estimate:
\bay
\label{p}
\|W\|_{\bS_1}\le\|\Psi\|_{L^\be\otimes_{\rm h}\!L^\be\otimes^{\rm h}\!L^\be}
\|T\|_{\bS_1}\|R\|,
\ey
where $\|\Psi\|_{L^\be\otimes_{\rm h}\!L^\be\otimes^{\rm h}\!L^\be}$ is the infimum of $\|\{\a_j\}_{j\ge0}\|_{L^\be(\ell^2)}\|\{\b_k\}_{k\ge0}\|_{L^\be(\ell^2)}
\|\{\g_{jk}\}_{j,k\ge0}\|_{L^\be(\B)}$
over all representations in \rf{yaH}.

\medskip

{\bf Definition 2.} {\it The Haagerup-like tensor product 
$L^\be(E_1)\otimes^{\rm h}\!L^\be(E_2)\otimes_{\rm h}\!L^\be(E_3)$ of the second kind
consists, be definition, of functions 
$\Psi$ that admit a representation
$$
\Psi(x_1,x_2,x_3)=\sum_{j,k\ge0}\a_{jk}(x_1)\b_{j}(x_2)\g_k(x_3)
$$
where $\{\b_j\}_{j\ge0},~\{\g_k\}_{k\ge0}\in L^\be(\ell^2)$, 
$\{\a_{jk}\}_{j,k\ge0}\in L^\be(\mB)$.}

\medskip

Suppose now that $T$ is a bounded linear operator, and $R\in\bS_1$. Then
the continuous linear functional 
$$
Q\mapsto
\trace\left(\left(
\iiint\Psi(x_1,x_2,x_3)\,dE_3(x_3)Q\,dE_1(x_1)T\,dE_2(x_2)
\right)R\right)
$$
on the class of compact operators determines a trace class operator,
which we call the triple operator integral
$$
W\df\upint\!\!\!\iint\Psi(x_1,x_2,x_3)\,dE_1(x_1)T\,dE_2(x_2)R\,dE_3(x_3).
$$
It follows easily from \rf{opno} that
\bay
\label{v}
\|W\|_{\bS_1}\le
\|\Psi\|_{L^\be\otimes^{\rm h}\!L^\be\otimes_{\rm h}\!L^\be}
\|T\|\cdot\|R\|_{\bS_1}.
\ey

Note that the above definitions of triple operator integrals extend the definition given in \cite{Pe6} in terms of the projective tensor product of the $L^\be$ spaces.

\

\section{Commutators of functions of almost commuting self-adjoint operators}
\label{comf}

\

In this section we prove that for a pair of almost commuting self-adjoint operators $A$ and $B$, the functional calculus
$$
\f\mapsto\f(A,B),\quad\f\in B_{\be,1}^1(\R^2),
$$
satisfies properties (i)--(iii) listed in the Introduction. Property (iv) will be established in the next section.

Properties (i) and (ii) are obvious.
To prove property (iii), we are going to use triple operator integrals with integrands in the Haagerup-like tensor products of $L^\be$ spaces and we find a representation of commutators $[\f(A,B),\psi(A,B)]$ in terms of such triple operator integrals.

Given a differentiable function $\f$ on $\R^2$, we define the divided differences
$\dg^{[1]}\f$ and $\dg^{[2]}\f$ on $\R^3$ by
$$
\big(\dg^{[1]}\f)(x_1,x_2,y\big)\df\frac{\f(x_1,y)-\f(x_2,y)}{x_1-x_2}\quad\mbox{and}\quad
\big(\dg^{[2]}\f\big)(x,y_1,y_2)=\frac{\f(x,y_1)-\f(x,y_2)}{y_1-y_2}.
$$

It was shown in \cite{ANP} and \cite{AP2} that if $\f$ is a bounded function on $\R^2$ whose Fourier transform is supported 
in the ball $\{\xi\in\R^2:~\|\xi\|\le1\}$, then
$$
\big(\dg^{[1]}\f)(x_1,x_2,y\big)=
\sum_{j,k\in\Z}\frac{\sin(x_1-j\pi)}{x_1-j\pi}\cdot\frac{\sin(x_2-k\pi)}{x_2-k\pi}
\cdot\frac{\f(j\pi,y)-\f(k\pi,y)}{j\pi-k\pi}.
$$
$$
\sum_{j\in\Z}\frac{\sin^2(x_1-j\pi)}{(x_1-j\pi)^2}
=\sum_{k\in\Z}\frac{\sin^2(x_2-k\pi)}{(x_2-k\pi)^2}=1,
\quad x_1~x_2\in\R,
$$
and
$$
\sup_{y\in\R}\left\|\left\{\frac{\f(j\pi,y)-\f(k\pi,y)}{j\pi-k\pi}
\right\}_{j,k\in\Z}\right\|_\B\le\const\|f\|_{L^\be(\R)}.
$$
Note that for $j=k$, we assume that 
$$
\frac{\f(j\pi,y)-\f(k\pi,y)}{j\pi-k\pi}
=\frac{\partial\f(x,y)}{\partial x}\Big|_{(j\pi,y)}.
$$
It follows that for such functions $\f$,
$$
\dg^{[1]}\f\in L^\be\otimes_{\rm h}\!L^\be\otimes^{\rm h}\!L^\be\quad
\mbox{and}\quad
\|\dg^{[1]}\f\|_{L^\be\otimes_{\rm h}\!L^\be\otimes^{\rm h}\!L^\be}
\le\const\|\f\|_{L^\be}.
$$
By rescaling, we can deduce from this that for bounded functions $\f$ on
$\R^2$ whose Fourier transform is supported in $\{\xi\in\R^2:~\|\xi\|\le\s\}$, we have
$$
\dg^{[1]}\f\in L^\be\otimes_{\rm h}\!L^\be\otimes^{\rm h}\!L^\be\quad
\mbox{and}\quad
\|\dg^{[1]}\f\|_{L^\be\otimes_{\rm h}\!L^\be\otimes^{\rm h}\!L^\be}
\le\const\s\|\f\|_{L^\be}.
$$

Suppose now that $\f\in B_{\be,1}^1(\R^2)$. Representing $\dg^{[1]}\f$ as
$$
\dg^{[1]}\f=\sum_{n\in\Z}\dg^{[1]}\f_n
$$
(see \S~\ref{intro}), we see that
$$
\dg^{[1]}\f\in L^\be\otimes_{\rm h}\!L^\be\otimes^{\rm h}\!L^\be\quad
\mbox{and}\quad
\|\f\|_{L^\be\otimes_{\rm h}\!L^\be\otimes^{\rm h}\!L^\be}
\le\const\|\f\|_{B_{\be,1}^1}.
$$

Similarly, for $\f\in B_{\be,1}^1(\R^2)$,
$$
\dg^{[2]}\f\in L^\be\otimes^{\rm h}\!L^\be\otimes_{\rm h}\!L^\be
\quad\mbox{and}\quad
\|\f\|_{L^\be\otimes^{\rm h}\!L^\be\otimes_{\rm h}\!L^\be}
\le\const\|\f\|_{B_{\be,1}^1}
$$
(see \cite{ANP} and \cite{AP2}).

\begin{thm}
\label{komut}
Let $A$ and $B$ be self-adjoint operators and let $Q$ be a bounded linear operator such that $[A,Q]\in\bS_1$ and $[B,Q]\in\bS_1$. 
Suppose that $\f\in B_{\be,1}^1(\R^2)$.
Then $[\f(A,B),Q\big]\in\bS_1$,
\begin{align}
\label{komQ}
\big[\f(A,B),Q\big]&=
\upint\!\!\!\iint\frac{\f(x,y_1)-\f(x,y_2)}{y_1-y_2}\,dE_A(x)\,dE_B(y_1)[B,Q]\,dE_B(y_2)\nonumber
\\[.2cm]
&+
\iint\!\!\upint\frac{\f(x_1,y)-\f(x_2,y)}{x_1-x_2}\,dE_A(x_1)[A,Q]\,dE_A(x_2)\,dE_B(y)
\end{align}
and
\bay
\label{ABQ}
\big\|[\f(A,B),Q\big]\big\|_{\bS_1}
\le\const\|\f\|_{B_{\be,1}^1(\R^2)}\big(\big\|[A,Q]\big\|_{\bS_1}+
\big\|[B,Q]\big\|_{\bS_1}\big).
\ey
\end{thm}

\Pf Let us first prove formula \rf{komQ} under the assumption that the divided differences $\dg_1\f$ and $\dg_2\f$ belong to the projective tensor products
$L^\be(E_A)\hat\otimes L^\be(E_A)\hat\otimes L^\be(E_B)$ and
$L^\be(E_A)\hat\otimes L^\be(E_B)\hat\otimes L^\be(E_B)$.
We have
\begin{align*}
\upint\!\!\!\iint&\frac{\f(x,y_1)-\f(x,y_2)}{y_1-y_2}\,dE_A(x)\,dE_B(y_1)(BQ-QB)\,dE_B(y_2)\\[.2cm]
&=\iiint\frac{\f(x,y_1)-\f(x,y_2)}{y_1-y_2}\,dE_A(x)\,dE_B(y_1)BQ\,dE_B(y_2)\\[.2cm]
&-\iiint\frac{\f(x,y_1)-\f(x,y_2)}{y_1-y_2}\,dE_A(x)\,dE_B(y_1)QB\,dE_B(y_2)\\[.2cm]
&=\iiint\frac{\f(x,y_1)-\f(x,y_2)}{y_1-y_2}y_1\,dE_A(x)\,dE_B(y_1)Q\,dE_B(y_2)\\[.2cm]
&-\iiint\frac{\f(x,y_1)-\f(x,y_2)}{y_1-y_2}y_2\,dE_A(x)\,dE_B(y_1)Q\,dE_B(y_2)\\[.2cm]
&=\iiint\big(\f(x,y_1)-\f(x,y_2)\big)\,dE_A(x)\,dE_B(y_1)Q\,dE_B(y_2)\\[.2cm]
&=\left(\iint\f(x,y_1)\,dE_A(x)\,dE_B(y_1)\right)Q
-\iint\f(x,y_2)\,dE_A(x)Q\,dE_B(y_2)\\[.2cm]
&=\f(A,B)Q-\iint\f(x,y)\,dE_A(x)Q\,dE_B(y).
\end{align*}
Similarly,
\begin{align*}
\iint\!\!\upint&\frac{\f(x_1,y)-\f(x_2,y)}{x_1-x_2}\,dE_A(x_1)[A,Q]\,dE_A(x_2)\,dE_B(y)\\[.2cm]
&=\iint\f(x,y)\,dE_A(x)Q\,dE_B(y)-Q\f(A,B)
\end{align*}
which proves \rf{komQ} under the above assumption.

Clearly, \rf{p}, \rf{v} and \rf{komQ} imply that under the above assumption
inequality \rf{ABQ} holds.

Suppose now that $\f$ is an arbitrary function in $B_{\be,1}^1(\R^2)$.
Representing $[\f(A,B),Q\big]$ in the form
$$
[\f(A,B),Q\big]=\sum_{n\in\Z}[\f_n(A,B),Q\big]
$$
(see \S~\ref{intro}), we find that it suffices to prove \rf{komQ} and \rf{ABQ}
for each function $\f_n$. 

As we have mentioned in Subsection 2.1, $\f_n$ is a restriction of an entire function of two variables to 
$\R\times\R$. Thus it suffices to prove \rf{komQ} and \rf{ABQ} in the case when $\f$ is an entire function. To complete the proof, we show that for entire functions $\f$ the divided differences $\dg^{[1]}\f$ and $\dg^{[1]}\f$ must belong to the projective tensor product
$L^\be\hat\otimes L^\be\hat\otimes L^\be$.  

Being an entire function, $\f$ admits an expansion
$$
\f(x,y)=\sum_{j=0}^\be\Big(\sum_{k=0}^\be a_{jk}x^jy^k\Big).
$$
Let $R$ be a positive number such that the spectra $\s(A)$ and $\s(B)$ are contained in $[-R/2,R/2]$.
Clearly,
$$
\|\f\|_{L^\be\hat\otimes L^\be}\le\sum_{j=0}^\be\Big(\sum_{k=0}^\be |a_{jk}|R^{j+k}\Big)
<\be
$$
and
\begin{align*}
\left\|\dg^{[1]}\f\right\|_{L^\be\hat\otimes L^\be\hat\otimes L^\be}&=
\left\|\sum_{j=0}^\be\left(\sum_{k=1}^\be
\Big(\sum_{l=0}^{j-1} a_{jk}x_1^lx_2^{j-1-l}y^k\Big)
\right)\right\|_{L^\be\hat\otimes L^\be\hat\otimes L^\be}\\[.2cm]
&\le\sum\limits_{j=0}^\be\left(\sum\limits_{k=1}^\be j|a_{jk}|R^{j+k-1}\right)<+\be,
\end{align*}
where in the above expressions $L^\be$ means $L^\be[-R,R]$. 

\medskip

The proof for $\dg^{[2]}\f$ is the same.
This completes the proof. $\bl$

To obtain the main result of the paper, we apply Theorem \ref{komut} in the case
$Q=\psi(A,B)$, where $\psi\in B^1_{\be,1}(\R^2)$.

\begin{thm}
\label{glav}
Let $A$ and $B$ be almost commuting self-adjoint operators and let $\f$ and $\psi$ be functions in the Besov class $B_{\be,1}^1(\R^2)$. Then
\begin{align}
\label{kom}
\!\!\!\!\!\big[\f(A,B),\psi(A,B)\big]\!&=
\upint\!\!\!\iint
\frac{\f(x,y_1)-\f(x,y_2)}{y_1-y_2}\,dE_A(x)\,dE_B(y_1)[B,\psi(A,B)]\,dE_B(y_2)\nonumber
\\[.2cm]
&+\!\!
\iint\!\!\upint\!
\frac{\f(x_1,y)-\f(x_2,y)}{x_1-x_2}dE_A(x_1)[A,\psi(A,B)]dE_A(x_2)dE_B(y)
\end{align}
and
\bay
\label{ner}
\big\|[\f(A,B),\psi(A,B)\big]\big\|_{\bS_1}
\le\const\|\f\|_{B_{\be,1}^1(\R^2)}\|\psi\|_{B_{\be,1}^1(\R^2)}
\big\|[A,B]\big\|_{\bS_1}.
\ey
\end{thm}

\medskip

{\bf Remark.}
Note that the right-hand side of inequality \rf{ner} does not involve the norms of $A$ or $B$. Suppose now that $A$ and $B$ are not necessarily bounded self-adjoint operators 
such that the operator $AB-BA$ is densely defined and extends to a trace class operator.
In this case formulae \rf{komQ} and \rf{kom} allow us to formally define commutators
$\big[\f(A,B),\psi(A,B)\big]$ by formula \rf{kom}. Formula \rf{kom} involves
$[B,\psi(A,B)]$ and $[A,\psi(A,B)]$ that can be formally defined by \rf{komQ}.
This allows us to give a definition of $\big[\f(A,B),\psi(A,B)\big]$ in the case when $A$ and $B$ do not have to be bounded. Inequality \rf{ner} still holds for such operators.

\

\section{An extension of the Helton--Howe trace formula}
\label{HelHow}

\

In this section we use the results of the previous section to extend the Helton--Howe trace formula.

\begin{thm}
Let $A$ and $B$ be almost commuting self-adjoint operators and let $\f$ and $\psi$ be functions in the Besov class $B_{\be,1}^1(\R^2)$. Then the following formula holds:
\begin{align}
\label{exHH}
\trace\big({\rm i}\big(\f(A,B)\psi(A,B)&-\psi(A,B)\f(A,B)\big)\big)\nonumber\\[.2cm]
&=\frac{1}{2\pi}
\iint_{\R^2}\left(\frac{\partial\f}{\partial x}\frac{\partial\psi}{\partial y}-
\frac{\partial\f}{\partial y}\frac{\partial\psi}{\partial x}\right)g(x,y)\,dx\,dy,
\end{align}
where $g$ is the Pincus principal function associated with the operators 
$A$ and $B$.
\end{thm}

\Pf It was proved in \cite{HH} that formula \rf{exHH} holds for infinitely differentiable functions. In particlar,
\begin{align*}
\trace\big({\rm i}\big(\f_n(A,B)\psi_m(A,B)&-\psi_m(A,B)\f_n(A,B)\big)\big)\\[.2cm]
&=\frac{1}{2\pi}
\iint_{\R^2}\left(\frac{\partial\f_n}{\partial x}\frac{\partial\psi_m}{\partial y}-
\frac{\partial\f_n}{\partial y}\frac{\partial\psi_m}{\partial x}\right)g(x,y)\,dx\,dy,
\end{align*}
where $\f_n=\f*W_n$ and $\psi_m=\psi*W_m$, see \S~\ref{Pre}. The results follows now from the obvious facts:
\begin{align*}
\sum_{m,n\in\Z}
\trace\big({\rm i}\big(\f_n(A,B)\psi_m(A,B)&-\psi_m(A,B)\f_n(A,B)\big)\big)\\[.2cm]
&=\trace\big({\rm i}\big(\f(A,B)\psi(A,B)-\psi(A,B)\f(A,B)\big)\big)
\end{align*}
and
\begin{align*}
\sum_{m,n\in\Z}
\iint_{\R^2}&\left(\frac{\partial\f_n}{\partial x}\frac{\partial\psi_m}{\partial y}-
\frac{\partial\f_n}{\partial y}\frac{\partial\psi_m}{\partial x}\right)g(x,y)\,dx\,dy
\\[.2cm]
&=
\iint_{\R^2}\left(\frac{\partial\f}{\partial x}\frac{\partial\psi}{\partial y}-
\frac{\partial\f}{\partial y}\frac{\partial\psi}{\partial x}\right)g(x,y)\,dx\,dy.
\quad\bl
\end{align*}
\medskip

It would be interesting to extend the notion of the Pincus principal function to the case of unbounded self-adjoint operators with trace class commutators and extend formula \rf{exHH} to unbounded almost commuting operators.

In the Introduction we associate with a pair $A$ and $B$ of almost commuting self-adjoint operator the operator $T$ defined by $A+{\rm i}B$ with trace class selfcommutator 
$[T^*,T]$. Then the operator $T$ is {\it essentially normal}, i.e., $[T^*,T]$ is compact. 
Consider its {\it essential spectrum} $\s_{\rm e}(T)$. It was proved in \cite{HH} that
on each component of $\C\setminus\s_{\rm e}(T)$ the principal function $g$ is constant and on a component of $\C\setminus\s_{\rm e}(T)$ it is equal to $-\ind(T-\l I)$, where $\l$ is a point in this component. This implies the following result:

\begin{thm}
\label{kompo}
Suppose that $A$ and $B$ are almost commuting self-adjoint operators
and $T=A+{\rm i}B$. Let $\O$ be a component of $\C\setminus\s_{\rm e}(T)$.
If $\f$ and $\psi$ are functions in $B_{\be,1}^1(\R^2)$ with supports in $\O$, then
\begin{align*}
\trace\big({\rm i}\big(\f(A,B)\psi(A,B)&-\psi(A,B)\f(A,B)\big)\big)\nonumber\\[.2cm]
&=-\frac{\ind(T-\l I)}{2\pi}
\iint_{\R^2}\left(\frac{\partial\f}{\partial x}\frac{\partial\psi}{\partial y}-
\frac{\partial\f}{\partial y}\frac{\partial\psi}{\partial x}\right)\,dx\,dy,
\end{align*}
where $\l$ is a point in $\O$.
\end{thm}

\

\section{Open problems}
\label{HelHow}

\

In this section we state two naturally arising problems.

\medskip

{\bf Problem 1.}
It was proved in \cite{Pe4} that for almost commuting self-adjoint operators $A$ and $B$, the functional calculus
$$
\f\mapsto\f(A,B),\quad
\f\in \big(L^\be(\R)\hat\otimes B_{\be,1}^1(\R)\big)\bigcap\big(B_{\be,1}^1(\R)\hat\otimes L^\be(\R)\big),
$$
is {\it almost multiplicative}, i.e.,
$$
(\f\psi)(A,B)-\f(A,B)\psi(A,B)\in\bS_1,\quad\f,\psi
\in \big(L^\be(\R)\hat\otimes B_{\be,1}^1(\R)\big)\bigcap\big(B_{\be,1}^1(\R)\hat\otimes L^\be(\R)\big).
$$
Here $\hat\otimes$ stands for projective tensor product.

The Besov class $B_{\be,1}^1(\R^2)$ does not form an algebra.
However, functions $\f(A,B)$ of $A$ and $B$ depend only on the restrictions of
$\f$ to the cartesian product $\s(A)\times\s(B)$ of the spectra of $A$ and $B$.
It is well-known that the class of restrictions of functions in $B_{\be,1}^1(\R^2)$ 
to a compact subset of $\R^2$ forms an algebra.

It would be interesting to find out whether the functional calculus
$\f\mapsto\f(A,B)$, $\f\in B_{\be,1}^1(\R^2)$, is also almost multiplicative.

\medskip

{\bf Problem 2.} It was proved in \cite{Pe4} that for almost commuting self-adjoint operators $A$ and $B$, the functional calculus
$$
\f\mapsto\f(A,B),\quad
\f\in B_{\be,1}^1(\R)\hat\otimes B_{\be,1}^1(\R),
$$
has the following property:
$$
\big(\f(A,B)\big)^*-\ov{\f}(A,B)\in\bS_1,\quad \f\in B_{\be,1}^1(\R)\hat\otimes B_{\be,1}^1(\R).
$$

We would like to pose the problem to find out whether the same property holds for arbitrary functions $\f$ in $B_{\be,1}^1(\R^2)$.

\

\noindent
\begin{tabular}{p{8cm}p{8cm}}
A.B. Aleksandrov  &  V.V. Peller \\
St.Petersburg Branch  & Department of Mathematics  \\
Steklov Institute of Mathematics  & Michigan State University\\
Fontanka 27   & East Lansing, Michigan 48824 \\
 191023 St-Petersburg  & USA\\
 Russia
\end{tabular}

\end{document}

%% file: classstartscr.tex
\newcommand{\lb}{\linebreak}

\renewcommand{\a}{\alpha}
\renewcommand{\b}{\beta}
\newcommand{\g}{\gamma}

\newcommand{\z}{\zeta}

\renewcommand{\l}{\lambda}

\newcommand{\s}{\sigma}

\newcommand{\f}{\varphi}
\renewcommand{\o}{\omega}

\newcommand{\D}{\Delta}
\renewcommand{\L}{\Lambda}

\renewcommand{\O}{\Omega}

\newcommand{\B}{{\mathscr B}}

\newcommand{\E}{{\mathscr E}}

\newcommand{\F}{{\mathscr F}}

\newcommand{\h}{{\mathscr H}}

\newcommand{\X}{{\mathscr X}}
\newcommand{\Y}{{\mathscr Y}}

\newcommand{\C}{{\Bbb C}}
\newcommand{\T}{{\Bbb T}}

\newcommand{\R}{{\Bbb R}}
\newcommand{\Z}{{\Bbb Z}}

\newcommand{\0}{{\boldsymbol{0}}}

\newcommand{\m}{{\boldsymbol m}}
\newcommand{\bS}{{\boldsymbol S}}

\newcommand{\rf}[1]{(\ref{#1})}

\newcommand{\df}{\stackrel{\mathrm{def}}{=}}

\newcommand{\supp}{\operatorname{supp}}

\newcommand{\trace}{\operatorname{trace}}

\newcommand{\const}{\operatorname{const}}

\newcommand{\eeq}{\end{equation}}
\newcommand{\beq}{\begin{equation}}
\newcommand{\bay}{\begin{eqnarray}}
\newcommand{\ba}{\begin{align*}}
\newcommand{\ea}{\end{align*}}
\newcommand{\ey}{\end{eqnarray}}
\newcommand{\bey}{\begin{eqnarray*}}
\newcommand{\eey}{\end{eqnarray*}}

\newcommand{\be}{\infty}

\newcommand{\bl}{\blacksquare}
\newcommand{\ess}{\operatorname{ess}}
\newcommand{\ind}{\operatorname{ind}}

\newcommand{\Pf}{{\bf Proof. }}
\newcommand{\im}{\operatorname{Im}}

\newcommand{\ov}{\overline}

\newtheorem{thm}{\hspace{\parindent}Theorem}[section]